\newcommand{\R}{I\!\!R}
\newcommand{\C}{\mathbb C}
\newcommand{\Gr}{\mathrm{gr}}
\newcommand{\Hcal}{\mathcal H}

\documentclass[oneside,draft,letterpaper,11pt]{amsart}

\usepackage{times}
\usepackage{amssymb}

\title[Complementary Lagrangians in Infinite Dimensional Symplectic Spaces]{Complementary Lagrangians in Infinite Dimensional Symplectic Hilbert Spaces}
\author[P.\ Piccione]{Paolo Piccione}
\author[D.\ Tausk]{Daniel V.\ Tausk}
\address{Departamento de Matem\'atica,\hfill\break\indent
Instituto de Matem\'atica e Estat\'\i stica\hfill\break\indent
Universidade de S\~ao Paulo, \hfill\break\indent Caixa Postal 66281,
CEP 05315--970, SP\hfill\break\indent Brazil}
\email{piccione@ime.usp.br, tausk@ime.usp.br}
\urladdr{http://www.ime.usp.br/\~{}piccione}
\urladdr{http://www.ime.usp.br/\~{}tausk}
\thanks{The authors are partially sponsored by CNPq.
The authors wish to thank Prof.\ Kenro Furutani for providing
instructive suggestions on the topic;
the material in this paper was developed after his observation
that Lagrangians can be characterized in terms of unbounded self-adjoint operators
(see Lemma~\ref{thm:lemFurutani}).}


\date{October 2nd, 2004}

\begin{document}

% Theorems and such

\theoremstyle{plain}\newtheorem*{teo}{Theorem}
\theoremstyle{plain}\newtheorem{prop}{Proposition}
\theoremstyle{plain}\newtheorem{lem}[prop]{Lemma}
\theoremstyle{plain}\newtheorem{cor}[prop]{Corollary}
\theoremstyle{definition}\newtheorem{defin}[prop]{Definition}
\theoremstyle{remark}\newtheorem{rem}[prop]{Remark}
\theoremstyle{plain} \newtheorem{assum}[prop]{Assumption}
\theoremstyle{definition}\newtheorem{example}{Example}

%%%%%

\begin{abstract}
Using the Spectral Theorem for unbounded self-adjoint operators we prove that any countable
family of Lagrangian subspaces of a symplectic Hilbert space admits a common complementary Lagrangian.
\end{abstract}

\maketitle
%\tableofcontents

%%%%%%%%%%%%%%%%
%%%%%%%%%%%%%%%%
\begin{section}{Introduction}\label{sec:intro}

A {\em symplectic Hilbert space\/} is a real Hilbert space $(V,\langle\cdot,\cdot\rangle)$ endowed with a symplectic form;
by a {\em symplectic form\/} we mean a bounded anti-symmetric bilinear form $\omega:V\times V\to\R$ that is represented by a
(anti-self-adjoint) linear isomorphism $H$ of $V$, i.e., $\omega=\langle H\cdot,\cdot\rangle$. If $H=PJ$ is the polar decomposition
of $H$ then $P$ is a positive isomorphism of $V$ and $J$ is a orthogonal complex structure on $V$;
the inner product $\langle P\cdot,\cdot\rangle$ on $V$ is therefore
equivalent to $\langle\cdot,\cdot\rangle$ and $\omega$ is represented by $J$ with respect to $\langle P\cdot,\cdot\rangle$.
We may therefore replace $\langle\cdot,\cdot\rangle$ with $\langle P\cdot,\cdot\rangle$ and assume from the beginning
that $\omega$ is represented by a orthogonal complex structure $J$ on $V$.
A subspace $S$ of $V$ is called {\em isotropic\/} if $\omega$ vanishes on $S$ or, equivalently, if $J(S)$ is contained in $S^\perp$.
A {\em Lagrangian subspace\/} of $V$ is a maximal isotropic subspace of $V$. We have that $L\subset V$ is Lagrangian
if and only if $J(L)=L^\perp$. If $L\subset V$ is Lagrangian then a Lagrangian $L'\subset V$ such that $V=L\oplus L'$
is called a {\em complementary Lagrangian\/} to $L$. Obviously every Lagrangian $L$ has a complementary Lagrangian,
namely, its orthogonal complement $L^\perp$. Given a pair $L_1$, $L_2$ of Lagrangians,
there are known sufficient conditions for the existence of a common complementary Lagrangian to $L_1$ and $L_2$ (see, for instance,
\cite{Fu}). In this paper we prove the following:

\begin{teo}
If $(V,\langle\cdot,\cdot\rangle,\omega)$ is a Hilbert symplectic
space then any countable family of Lagrangian subspaces of\/ $V$
admits a common complementary Lagrangian.
\end{teo}

The existence of a common complementary Lagrangian is proven first in the case
of two Lagrangians $L$ and $L_1$ such that $L\cap L_1=\{0\}$ (Corollary~\ref{thm:casointerzero}).
In this case $L$ is the graph of a densely defined self-adjoint operator on
$L_1^\perp$ (Lemma~\ref{thm:lemFurutani}), and the result is obtained as an
application of the spectral theorem (Lemma~\ref{thm:AAlinha} and Lemma~\ref{thm:outroAAlinha}).
The existence of a common complementary Lagrangians is then proven
in the general case  by a reduction argument (Proposition~\ref{thm:commoncompl}), and the final result
is an application of Baire's category theorem.

\end{section}

\begin{section}{Proof of the Result}
\label{sec:proof}

In what follows, $(V,\langle\cdot,\cdot\rangle,\omega)$ will denote a Hilbert symplectic space such that
$\omega$ is represented by a orthogonal complex structure $J$ on $V$. We will denote by $\Lambda(V)$ the set
of all Lagrangian subspaces of $V$. It follows from Zorn's Lemma that $V$ indeed has Lagrangian subspaces, i.e.,
$\Lambda(V)\ne\emptyset$. Given $L_0,L_1\in\Lambda(V)$ then $(L_0+L_1)^\perp=J(L_0\cap L_1)$; in particular,
$L_0\cap L_1=\{0\}$ if and only if $L_0+L_1$ is dense in $V$. For $L\in\Lambda(V)$, we denote by $\mathcal O(L)$
the subset of $\Lambda(V)$ consisting of Lagrangians complementary to $L$.
Given a real Hilbert space $\Hcal$, we denote by $\Hcal^\C$ the orthogonal direct sum $\Hcal\oplus\Hcal$
endowed with the orthogonal complex structure $J$ defined by $J(x,y)=(-y,x)$. If $A:D\subset\Hcal\to\Hcal$ is a densely defined linear
operator on $\Hcal$ then $J\big(\Gr(A)^\perp\big)=\Gr(A^*)$. It follows that $\Gr(A)$ is Lagrangian in $\Hcal^\C$
if and only if $A$ is self-adjoint; in this case, $\Gr(A)$ is complementary to $\{0\}\oplus\Hcal$ if and only
if $A$ is bounded.
\begin{lem}\label{thm:lemFurutani}
Given $L\in\Lambda(\Hcal^\C)$ with $L\cap\big(\{0\}\oplus\Hcal\big)=\{0\}$ then $L$ is the graph of a densely defined
self-adjoint operator $A:D\subset\Hcal\to\Hcal$.
\end{lem}
\begin{proof}
The sum $L+\big(\{0\}\oplus\Hcal\big)$ is dense in $\Hcal^\C$; thus, denoting by $\pi_1:\Hcal^\C\to\Hcal$ the projection
onto the first summand, we have that $D=\pi_1(L)=\pi_1\big(L+\big(\{0\}\oplus\Hcal\big)\big)$
is dense in $\Hcal$. Hence $L$
is the graph of a densely defined operator $A:D\to\Hcal$, which is self-adjoint by the remarks above.
\end{proof}

Given Lagrangians $L_0,L_1\in\Lambda(V)$ with $V=L_0\oplus L_1$ then we have an isomorphism $\rho_{L_1,L_0}:L_1\to L_0$
defined by $\rho_{L_1,L_0}=P_{L_0}\circ J\vert_{L_1}$, where $P_{L_0}$ denotes the orthogonal projection onto $L_0$.
The map:
\begin{equation}\label{eq:rhosymplecto}
V=L_0\oplus L_1\ni x+y\longmapsto\big(x,-\rho_{L_1,L_0}(y)\big)\in L_0\oplus L_0=L_0^\C
\end{equation}
is a symplectomorphism, i.e., it is
an isomorphism that preserves the symplectic forms. Thus, we get a one-to-one correspondence $\varphi_{L_0,L_1}$
between Lagrangian subspaces $L$ of $V$ with $L\cap L_1=\{0\}$ and densely defined self-adjoint operators $A:D\subset L_0\to L_0$;
more explicitly, we set $A=\varphi_{L_0,L_1}(L)$ if the map \eqref{eq:rhosymplecto} carries
$L$ to the graph of $-A$.

\begin{lem}\label{thm:AAlinha}
Let $L_0,L_1,L,L'\in\Lambda(V)$ be Lagrangians such that $L_0$ and $L'$ are complementary to $L_1$ and $L\cap L_1=\{0\}$.
Set $\varphi_{L_0,L_1}(L)=A:D\subset L_0\to L_0$ and $\varphi_{L_0,L_1}(L')=A':L_0\to L_0$. Then
$L'$ is complementary to $L$ if and only if $(A-A'):D\to L_0$ is an isomorphism.
\end{lem}
\begin{proof}
The map \eqref{eq:rhosymplecto} carries $L$ and $L'$ respectively to $\Gr(-A)$ and $\Gr(-A')$. We thus
have to show that $L_0^\C=\Gr(-A)\oplus\Gr(-A')$ if and only if $A-A'$ is an isomorphism. This follows by observing
that $(x,y)=(u,-Au)+(u',-A'u')$ is equivalent to $\big(u+u',(A'-A)u\big)=(x,y+A'x)$, for all $x,y,u'\in L_0$, $u\in D$.
\end{proof}

\begin{lem}\label{thm:outroAAlinha}
If $A:D\subset\Hcal\to\Hcal$ is a densely defined self-adjoint operator then for every $\varepsilon>0$ there exists a bounded
self-adjoint operator $A':\Hcal\to\Hcal$ with $\Vert A'\Vert\le\varepsilon$ and such that $(A-A'):D\to\Hcal$ is an isomorphism.
\end{lem}
\begin{proof}
By the Spectral Theorem for unbounded self-adjoint operators, we may assume that
$\Hcal=L^2(X,\mu)$ and $A=M_f$, where $(X,\mu)$ is a measure space, $f:X\to\R$ is a measurable function
and $M_f$ denotes the multiplication operator by $f$ defined on $D=\big\{\phi\in L^2(X,\mu):f\phi\in L^2(X,\mu)\big\}$.
In this situation, the operator $A'$ can be defined as $A'=M_g$, where $g=\varepsilon\cdot\chi_\varepsilon$
and $\chi_\varepsilon$ is
the characteristic function of the set $f^{-1}\big([-\frac\varepsilon2,\frac\varepsilon2]\big)$;
clearly $\Vert A'\Vert\le\Vert g\Vert_\infty=\varepsilon$.
The conclusion follows by observing that $A-A'=M_{f-g}$, and $\vert f-g\vert\ge\frac\varepsilon2$ on $X$.
\end{proof}

\begin{cor}\label{thm:casointerzero}
Given $L_1,L\in\Lambda(V)$ with $L_1\cap L=\{0\}$ then there exists a common complementary Lagrangian $L'\in\Lambda(V)$
to $L_1$ and $L$.
\end{cor}
\begin{proof}
Set $L_0=L_1^\perp$ and $A=\varphi_{L_0,L_1}(L)$. Lemma~\ref{thm:outroAAlinha} gives us a bounded self-adjoint operator
$A':L_0\to L_0$ with $A-A'$ an isomorphism. Set $L'=\varphi_{L_0,L_1}^{-1}(A')$; $L'$ is a Lagrangian complementary
to $L_1$, because $A'$ is bounded. It is also complementary to $L$, by Lemma~\ref{thm:AAlinha}.
\end{proof}

If $V=V_1\oplus V_2$ is a orthogonal direct sum decomposition into $J$-invariant subspaces $V_1$ and $V_2$, then
$V_1$ and $V_2$ are symplectic Hilbert subspaces of $V$. Given subspaces $L_1\subset V_1$ and $L_2\subset V_2$ then
$L_1\oplus L_2$ is Lagrangian in $V$ if and only if $L_i$ is Lagrangian in $V_i$, for $i=1,2$. A Lagrangian
subspace $L\in\Lambda(V)$ is of the form $L=L_1\oplus L_2$ with $L_i\in\Lambda(V_i)$, $i=1,2$, if and only if
$L$ is invariant by the orthogonal projection $P_{V_1}$ onto $V_1$. In this case, $L_i=P_{V_i}(L)=L\cap V_i$, $i=1,2$.
If $S$ is a closed isotropic subspace of $V$ then a decomposition $V=V_1\oplus V_2$ of the type above can be obtained
by setting $V_1=S\oplus J(S)$ and $V_2=V_1^\perp$. Then, if $L\in\Lambda(V)$ contains $S$, it follows that
$P_{V_1}(L)=S$; namely, $S\subset L$ implies $L\subset J(S)^\perp$ and $J(S)^\perp$ is invariant
by $P_{V_1}$. Hence $L=S\oplus P_{V_2}(L)$.
\begin{prop}\label{thm:commoncompl}
Given $L,L'\in\Lambda(V)$ then $\mathcal O(L)\cap\mathcal O(L')\ne\emptyset$.
\end{prop}
\begin{proof}
Set $S=L\cap L'$, $V_1=S\oplus J(S)$, and $V_2=V_1^\perp$. Then $L=S\oplus P_{V_2}(L)$, $L'=S\oplus P_{V_2}(L')$,
and $P_{V_2}(L)\cap P_{V_2}(L')=(L\cap V_2)\cap(L'\cap V_2)=\{0\}$. By Corollary~\ref{thm:casointerzero},
there exists a Lagrangian $R\in\Lambda(V_2)$ complementary to both $P_{V_2}(L)$ and $P_{V_2}(L')$ in $V_2$.
Hence $J(S)\oplus R\in\Lambda(V)$ is in $\mathcal O(L)\cap\mathcal O(L')$.
\end{proof}

The map $L\mapsto P_L$ is a bijection from $\Lambda(V)$ onto the space of bounded self-adjoint maps $P:V\to V$
with $P^2=P$ and $PJ+JP=J$. Such bijection induces a topology on $\Lambda(V)$ which makes it homeomorphic to a complete
metric space. Moreover, for any $L_0,L_1\in\Lambda(V)$ with $V=L_0\oplus L_1$, the set $\mathcal O(L_1)$ is open
in $\Lambda(V)$ and the map $\mathcal O(L_1)\ni L\mapsto\varphi_{L_0,L_1}(L)$ is a homeomorphism onto
the space of bounded self-adjoint operators on $L_0$.

\begin{lem}
For any $L_0\in\Lambda(V)$, the set $\mathcal O(L_0)$ is dense in $\Lambda(V)$.
\end{lem}
\begin{proof}
Given $L\in\Lambda(V)$, Proposition~\ref{thm:commoncompl} gives us $L_1\in\mathcal O(L_0)\cap\mathcal O(L)$. By Lemma~\ref{thm:outroAAlinha},
the bounded self-adjoint operator $A=\varphi_{L_0,L_1}(L)$ on $L_0$ is the limit of a sequence of bounded self-adjoint isomorphisms
$A_n:L_0\to L_0$. Hence the sequence $\varphi_{L_0,L_1}^{-1}(A_n)$ is in $\mathcal O(L_0)$ and it tends to $L$.
\end{proof}

\begin{proof}[Proof of Theorem]
Let $(L_n)_{n\ge1}$ be a sequence in $\Lambda(V)$. Each $\mathcal O(L_n)$ is open and dense in $\Lambda(V)$,
hence $\bigcap_{n=1}^\infty\mathcal O(L_n)$ is dense in $\Lambda(V)$, by Baire's category theorem.
\end{proof}

\end{section}

\end{document}